\documentclass[11pt]{amsart}

\usepackage{amssymb,geometry,stmaryrd,color}
\usepackage{amsmath}
\usepackage{amsfonts}
\usepackage{fullpage}

\usepackage{hyperref}
\hypersetup{colorlinks,linkcolor={red},citecolor={blue},urlcolor={blue}}

\usepackage{graphicx}
\usepackage{caption}
\begin{document}

\newtheorem{theorem}{Theorem}[section]
\newtheorem{lemma}[theorem]{Lemma}
\newtheorem{proposition}[theorem]{Proposition}
\newtheorem{corollary}[theorem]{Corollary}
\newtheorem{conjecture}[theorem]{Conjecture}
\newtheorem{example}{Example}

\newtheorem{definition}[theorem]{Definition}
\newtheorem{remark}[theorem]{Remark}
\newtheorem{notation}[theorem]{Notation}
\newtheorem{question}[theorem]{Question}

\numberwithin{equation}{section}

\def\s{{\bf s}} 
\def\t{{\bf t}} 
\def\u{{\bf u}} 
\def\x{{\bf x}} 
\def\y{{\bf y}} 
\def\z{{\bf z}} 
\def\B{{\bf B}} 
\def\C{{\bf C}} 
\def\D{{\bf D}}
\def\K{{\bf K}}
\def\F{{\bf F}}
\def\M{{\bf M}}
\def\ML{{\bf ML}}
\def\Nn{{\bf N}}
\def\G{{\bf \Gamma}} 
\def\W{{\bf W}}
\def\X{{\bf X}}
\def\U{{\bf U}}
\def\V{{\bf V}}
\def\Un{{\bf 1}}
\def\Y{{\bf Y}}
\def\Z{{\bf Z}}
\def\P{{\bf P}}
\def\Q{{\bf Q}}
\def\S{{\bf S}}
\def\L{{\bf L}}
\def\T{{\bf T}}

\def\cB{{\mathcal{B}}} 
\def\cC{{\mathcal{C}}} 
\def\cD{{\mathcal{D}}} 
\def\cG{{\mathcal{G}}} 
\def\cK{{\mathcal{K}}} 
\def\cL{{\mathcal{L}}} 
\def\cR{{\mathcal{R}}} 
\def\cS{{\mathcal{S}}}
\def\cU{{\mathcal{U}}}
\def\cV{{\mathcal{V}}} 
\def\cX{{\mathcal X}}
\def\cY{{\mathcal Y}}
\def\cZ{{\mathcal Z}}

\def\Ea{E_\a}
\def\eps{{\varepsilon}} 
\def\esp{{\mathbb{E}}} 
\def\Ga{{\Gamma}}

\def\lacc{\left\{}
\def\lcr{\left[}
\def\lpa{\left(}
\def\lva{\left|}
\def\racc{\right\}}
\def\rpa{\right)}
\def\rcr{\right]}
\def\rva{\right|}

\def\prst{{\leq_{st}}}
\def\prost{{\prec_{st}}}
\def\prcvx{{\prec_{cx}}}
\def\Rr{{\bf R}}

\def\CC{{\mathbb{C}}}
\def\EE{{\mathbb{E}}}
\def\NN{{\mathbb{N}}} 
\def\QQ{{\mathbb{Q}}} 
\def\PP{{\mathbb{P}}}
\def\ZZ{{\mathbb{Z}}}
\def\RR{{\mathbb{R}}}

\def\Tt{{\bf \Theta}}
\def\Ttt{{\tilde \Tt}}

\def\a{\alpha}
\def\A{{\bf A}}
\def\AA{{\mathcal A}}
\def\hAA{{\hat \AA}}
\def\hL{{\hat L}}
\def\hT{{\hat T}}

\def\claw{\stackrel{(d)}{\longrightarrow}}
\def\elaw{\stackrel{(d)}{=}}
\def\pslaw{\stackrel{a.s.}{\longrightarrow}}
\def\qed{\hfill$\square$}

\newcommand*\pFqskip{8mu}
\catcode`,\active
\newcommand*\pFq{\begingroup
        \catcode`\,\active
        \def ,{\mskip\pFqskip\relax}%
        \dopFq
}
\catcode`\,12
\def\dopFq#1#2#3#4#5{%
        {}_{#1}F_{#2}\biggl[\genfrac..{0pt}{}{#3}{#4};#5\biggr]%
        \endgroup
}

\def\ii{{\rm i}}

\title[Non-SD of alpha Cauchy]{Self-decomposability of $\alpha$-Cauchy distributions}

\author[M.~Wang]{Min Wang}

\address{School of Mathematics and Statistics, Wuhan University of Technology,  Wuhan, 430063, China}

\email{minwangmath@whut.edu.cn}

\author[S.~Yin]{Sheng Yin}

\address{Institute for Advanced Study in Mathematics, Harbin Institute of Technology, Harbin, 150001, China}

\email{sheng.yin@hit.edu.cn}



\begin{abstract} 
In 2009, Yano, Yano and Yor asked whether the $\alpha$-Cauchy
distributions are self-decomposable or at least infinitely divisible in their study of hitting
times for symmetric stable processes. 
Recently, the first named author proved that $\alpha$-Cauchy
distribution is infinitely divisible if and only if $1 < \alpha \leq 2$. In this paper, we prove that $\alpha$-Cauchy
distribution is self-decomposable if and only if $\alpha = 2$. 
 The proof is based on a new criterion for self-decomposability in the symmetric case. As an application, we show that the first hitting time of a nonzero point by a symmetric stable process of index
$1<\alpha<2$ starting from zero is not self-decomposable. 
\end{abstract}

\maketitle

\section{Introduction}
A random variable $X$ is said to be self-decomposable if for every
$c\in(0,1)$, it can be written as
$$ X \elaw cX'+X_c,$$
where $X'$ is a copy of $X$ and is independent of $X_c$.
Self-decomposable
distributions form an important subclass of infinitely divisible
distributions. They can be viewed as limit distributions of
a class of Markov processes called processes of Ornstein–Uhlenbeck type. 
Both generalized gamma convolutions (see \cite{Bon92}) and stable distributions are self-decomposable. 
Well-known examples of self-decomposable distributions are the normal, Cauchy, gamma, exponential, Gumbel; see \cite{Sat99}. Further examples include the half-Cauchy distribution \cite{Die98} and the Fréchet distributions \cite{BS13}. Related questions concerning inverse beta distributions and weak variance generalized gamma convolutions are studied in \cite{BS15, BLM20}.
We refer to Sato \cite{Sat99} for their general properties.

In 2009, Yano, Yano and Yor \cite{YYY09} introduced the $\alpha$-Cauchy variable to study the first hitting times of points for one-dimensional symmetric stable L\'evy processes.  
The density of $\alpha$-Cauchy variable, denoted by $\mathcal{C}_\a$, is 
\begin{equation}\label{eq density alpha Cauchy}
    f_{\mathcal{C}_\a}(x) = \frac{\sin(\pi/\a)}{2\pi/\a}\frac{1}{1+|x|^\a}, \quad \alpha > 1, \,\, x \in \mathbb{R}. 
\end{equation}
Yano, Yano and Yor \cite[Remark 2.9(i)]{YYY09} posed the following problem.
\vspace{2mm}
\begin{center}
  \textbf{Problem:} \textit{Is it true that the law of $\mathcal{C}_\a$ is self-decomposable (or infinitely divisible at least)?}
\end{center}
\vspace{2mm}
The first named author \cite{Wan26} proved that $\mathcal{C}_\alpha$ is infinitely divisible if and only if $1 < \alpha \leq 2$. 
In this paper, we further study the self-decomposability of $\mathcal{C}_\a$. 

\begin{theorem}\label{Theorem non-SD of alpha Cauchy}
The $\alpha$-Cauchy variable $\mathcal{C}_\a$ is self-decomposable if and only if  $\alpha = 2$. 
\end{theorem}

The  $\alpha$-Cauchy variable and the first hitting times of points for one-dimensional symmetric stable L\'evy processes are closely related; cf. \cite{YYY09, LS14, LS19, Wan26}. 
Let $X_\alpha = (X_\alpha(t): t\geq 0)$ be the symmetric stable L\'evy process of index $\alpha$ starting from zero. 
Let $T_{\{a\}}(X_\alpha)$ denote the first hitting time of point $a \neq 0$ for $X_\alpha$. These hitting times
are almost surely finite when $1<\alpha\le2$; see~\cite{YYY09, LS14}.
As an application of Theorem \ref{Theorem non-SD of alpha Cauchy}, we establish the following result.

\begin{corollary}\label{Corollary}
    For $1<\alpha <2$, the first hitting time of point $a \neq 0$ for a symmetric stable L\'evy process of index $\alpha$ starting from zero, i.e. $T_{\{a\}}(X_\alpha)$, is not self-decomposable. 
\end{corollary}

When $\alpha = 2$, 
it’s already known that the first hitting times of nonzero points for one-dimensional Brownian motion have positive $1/2-$stable distributions, see formula (5.13) in \cite{YYY09}. Therefore, these laws are self-decomposable. 

In order to prove Theorem \ref{Theorem non-SD of alpha Cauchy}, we establish a new criterion for symmetric self-decomposable distributions.

\begin{theorem}\label{theorem new criterion SD}
Let $X$ be a symmetric random variable whose characteristic function
$\phi$ is strictly positive on $\mathbb{R}$. Let $\Psi:=-\log\phi$, and
assume that $\Psi\in C^1(\mathbb{R}\setminus\{0\})$ and
 $$\lim_{z\to0}z\Psi'(z)=0,
 \qquad \lim_{z\to+\infty}z\Psi'(z)=L\in(0,\infty).$$
Define 
\begin{equation}
    \label{eq def J}
    J(0) = 1, \quad \text{and} \quad J(z) := 1 - \frac{z \Psi'(z)}{L},\quad \text{for} \quad z \neq 0 .
\end{equation}
\begin{itemize}
    \item[(i)] $X$ is self-decomposable if and only if $J$ is a characteristic
function of a symmetric
probability measure with no atom at zero. 
    \item[(ii)] If $X$ is self-decomposable and  $\int_\mathbb{R}|J(z)|dz < \infty$, then $  \int_\mathbb{R} J(z)dz \geq 0$.
\end{itemize}
\end{theorem}

Note that the characteristic function of an infinitely divisible distribution never vanishes; see e.g. Lemma 7.5 in \cite{Sat99}. Consequently, a symmetric infinitely divisible distribution has a strictly positive characteristic function. Hence, the condition that $\phi$ is strictly positive on $\mathbb{R}$ in Theorem \ref{theorem new criterion SD} always holds for all symmetric self-decomposable distributions.

For $1< \alpha < 2$, we construct the function $J_\a$ corresponding to $\mathcal{C}_\alpha$ according to the steps in Theorem \ref{theorem new criterion SD}. We then prove that $J_\a$ is absolutely integrable and that its integral over $\mathbb{R}$ is
strictly negative, hence $\mathcal{C}_\alpha$ is not self-decomposable by the necessary part of Theorem \ref{theorem new criterion SD}. The proofs of these facts depend on the asymptotic properties of $\phi_\a(z)$ and $\phi'_\a(z)$, which are derived from a representation of the characteristic function of $\mathcal{C}_\a$ established in our previous work
 \cite[Lemma 2.7]{Wan26}.

\begin{lemma}
    \label{lemma asymptotic estimates alpha Cauchy}
    For $1<\alpha <2$, the characteristic function of the $\alpha$-Cauchy variable admits the representation
\begin{equation}
    \label{eq Laplace form for CF}
    \phi_\a(z) = k_\a  \int_0^\infty  e^{-|z|y}  \frac{ y^\alpha}{  y^{2\alpha} + 2\cos(\pi \alpha/2)y^\alpha + 1}dy, 
\end{equation}
where $k_\a$ is a normalization constant such that $ \phi_\a(0) = 1$.  Moreover, $\phi_\a(z)$ and $\phi'_\a(z)$ have the following asymptotic properties:
\begin{equation}
    \label{eq phi infty}
        \phi_\a (z) = k_\alpha \left[ \Gamma(1+\alpha)z^{-1-\alpha } - 2 \cos(\pi \alpha/2) \Gamma(1 + 2\alpha) z^{-1-2\alpha}\right]  + O(z^{-1-3\alpha}) , \quad z \rightarrow +\infty;
\end{equation}
\begin{align}
\label{eq phi' infty}
 \phi_\alpha'(z)
 =-k_\alpha
 \left[\Gamma(2+\alpha)z^{-2-\alpha}-2 \cos(\pi\alpha/2)\Gamma(2+2\alpha) z^{-2-2\alpha}\right]
       +O(z^{-2-3\alpha}), \quad z \rightarrow +\infty;
\end{align}
and 
\begin{equation}
\label{eq phi' sym at 0}
     \phi_\a'(z) \sim  - k_\a \Gamma(2-\alpha) z^{\alpha -2} ,  \quad z \downarrow 0.
\end{equation}
\end{lemma}

The rest of this paper is organized as follows. In Section \ref{sec proof criterion}, we prove Theorem \ref{theorem new criterion SD} by the relation between self-decomposable distributions and the Ornstein-Uhlenbeck processes. 
In Section  \ref{sec proof non-SD}, we give the details of the proof of our main Theorem \ref{Theorem non-SD of alpha Cauchy}. The asymptotic estimates in Lemma \ref{lemma asymptotic estimates alpha Cauchy} are proved
in Section \ref{sec proof lemma}, and the hitting time application in Corollary \ref{Corollary} is proved in Section \ref{sec proof corollary}.

\section{Proof of Theorem \ref{theorem new criterion SD}}
\label{sec proof criterion}
\textbf{Necessity in (i):}
Wolfe \cite{Wol82} and Jurek and Vervaat \cite{JV83} showed that the distribution of a random variable $X$ is self-decomposable if and only if
$$X \elaw \int_0^\infty e^{-s} \, dY_s,$$ 
for some L\'evy process $Y=(Y_s, s\ge 0)$ with $\mathbb{E}[\log(1\vee |Y_s|)]<\infty$ for all $s$.
By Theorem 17.5 in Sato \cite{Sat99}, we have 
\begin{equation}
    \Psi(z) =  \int_0^\infty \psi_{Y_1}(e^{-s}z)ds = \int_0^1 \psi_{Y_1}(uz) \frac{du}{u} = \int_0^z \psi_{Y_1}(t) \frac{dt}{t}, \quad z \in \mathbb{R},
\end{equation}
where $\Psi(z) = -\log \mathbb{E}[e^{izX}]$ and $\psi_{Y_1}(z) = - \log \mathbb{E}[e^{izY_1}]$. Then we have 
\begin{equation}
   z \Psi'(z) =   \psi_{Y_1}(z), \quad z \in \mathbb{R}.
\end{equation}
Because $X$ is symmetric, we have $\Psi(z) = \Psi(-z)$ and then $\psi_{Y_1}(z) = z \Psi'(z) = -z \Psi'(-z) =\psi_{Y_1}(-z)$. This implies that $\mathbb{E}[e^{izY_1}] = \mathbb{E}[e^{iz(-Y_1)}]$, thus $Y_1$ is a symmetric infinitely divisible random variable. By Theorem 8.1 and Exercise 18.1 in Sato \cite{Sat99}, we have 
\begin{equation}
    \label{LK of Y}
     z\Psi'(z) = \psi_{Y_1}(z) = \frac{1}{2}Az^2 + \int_{-\infty}^{\infty} (1-\cos (zx))\nu(dx),
\end{equation}
where $A \geq 0$ and $\nu$ is a symmetric measure on $\mathbb{R}$ satisfying 
$$\nu(\{0\}) = 0 \quad \text{and} \quad \int_{\mathbb{R}} \min (|x|^2, 1) \, \nu(dx) < \infty. $$
If $z\Psi'(z)$ is bounded, then $A = 0$. We next prove that $\nu(\mathbb{R})$ is finite. 
\begin{align}\label{eq DC1}
    \nu(\mathbb{R}) & = \int_\mathbb{R}  \lim_{M \rightarrow +\infty} (1 - \frac{\sin(Mx)}{Mx})  \nu(dx)  \\
    \label{eq DC2}
    &\leq \lim_{M \rightarrow +\infty}  \int_\mathbb{R} (1 - \frac{\sin(Mx)}{Mx})  \nu(dx)  \\
    & = \lim_{M \rightarrow +\infty}  \int_\mathbb{R} \left[1 - \frac{1}{M} \int_0^M\cos(zx)dz\right]  \nu(dx) \\
    & = \lim_{M \rightarrow +\infty}  \int_\mathbb{R} \left[ \frac{1}{M} \int_0^M (1-\cos(zx))dz\right]  \nu(dx)  \\
    \label{eq2}
    &= \lim_{M \rightarrow +\infty} \frac{1}{M} \int_0^M \left[  \int_\mathbb{R} (1-\cos(zx)) \nu(dx) \right] dz  \\
    &= \lim_{M \rightarrow +\infty} \frac{1}{M} \int_0^M  z \Psi'(z) dz \\
    \label{eq3}
    & = \lim_{z \rightarrow  \infty} z\Psi'(z) = L.
\end{align}
\eqref{eq DC2} follows from Fatou's lemma. \eqref{eq2} follows from Tonelli's theorem. \eqref{eq3} is guaranteed by the assumption that $\lim_{z \rightarrow  \infty} z\Psi'(z)$ exists. 

Now $\nu(\mathbb{R})$ is finite, we can return to \eqref{eq DC1} and use dominated convergence to obtain that 
\begin{equation}
    \nu(\mathbb{R})  = \int_\mathbb{R}  \lim_{M \rightarrow +\infty} (1 - \frac{\sin(Mx)}{Mx})  \nu(dx) = 
     \lim_{M \rightarrow +\infty}  \int_\mathbb{R} (1 - \frac{\sin(Mx)}{Mx})  \nu(dx) = L.
\end{equation}

Then the function 
\begin{align*}
      J(z) &:= 1 - \frac{z \Psi'(z)}{L} = 1 - \frac{1}{L}  \int_{-\infty}^{\infty} (1-\cos (zx))\nu(dx) \\
      &= \frac{1}{L}  \int_{-\infty}^{\infty} \cos (zx)\nu(dx) = \frac{1}{L}  \int_{-\infty}^{\infty} e^{izx}\nu(dx)
\end{align*}
Recall that $\nu$ is the L\'evy measure of $Y_1$ and it is symmetric, see \eqref{LK of Y}. We can conclude that $J(z)$ is the characteristic function of a symmetric probability law with no atom at zero.

\textbf{Sufficiency in (i):}
Conversely, suppose that $J$ is a characteristic
function of a symmetric
probability measure. For
$0<c<1$, direct computation gives
\begin{equation}\label{eq SD-ratio}
 \frac{\phi(z)}{\phi(cz)}
 =\exp\left\{L\int_c^1\frac{J(uz)-1}{u}d u\right\}.
\end{equation}
We want to prove that $\frac{\phi(z)}{\phi(cz)}$ is a characteristic
function. 
Let $U_c, \,0 < c < 1,\,$ be a random variable with density 
$\frac{\mathbf{1}_{(c,1)}(u)}{u\log(1/c)}\,$ and let $W$ be a random variable with characteristic
function $J$. The random variables $U_c$ and $W$ are independent. The characteristic
function of $U_c W$ is 
\begin{align}
    \mathbb{E}(e^{izU_c W}) &= \mathbb{E}[\mathbb{E}(e^{izU_c W}|U_c)] \\
    &= \mathbb{E} [ J(zU_c)] \\
    &=\frac{1}{\log(1/c)} \int_c^1 J(z u) \, \frac{du}{u}. 
\end{align}
Hence 
\begin{equation}
    \frac{\phi(z)}{\phi(cz)}
 =\exp\left\{L \log(1/c)\int_\mathbb{R} (e^{izx}-1) \nu_c(d x)\right\},
\end{equation}
 where $\nu_c$ is the distribution of $U_c W$. By L\'evy-khintchine representation of infinitely divisible distributions (see Theorem 8.1 in Sato \cite{Sat99}), $\frac{\phi(z)}{\phi(cz)}$ is a characteristic
function of an infinitely divisible distribution. 
 This proves self-decomposability.

\textbf{(ii)} By Theorem \ref{theorem new criterion SD} (i), if $X$ is self-decomposable, then $J$ is a characteristic function. If further
$\int_\mathbb{R}|J(z)|dz < \infty$, then $\nu/L$ has a bounded continuous density 
$$f(y) = \frac{1}{2\pi} \int_\mathbb{R} e^{-izy} J(z)dz,$$
see e.g. Durrett \cite[Theorem 3.3.14]{Dur19} for its proof. In particular, setting $y = 0$, we have 
$$  \int_\mathbb{R} J(z)dz  =  2\pi f(0) \geq 0. $$
This completes the proof. 

\begin{remark}
    In Theorem \ref{theorem new criterion SD} (ii), if X is symmetric and self-decomposable, the function $J$ defined as \eqref{eq def J} is not automatically in $L^1(\mathbb{R})$. For example, the characteristic function of the Linnik variable $\Lambda_\a$ of index $0<\alpha \leq 2$ is 
    \begin{equation}
        \mathbb{E}[e^{i\theta\Lambda_\a}] = \frac{1}{1+|\theta|^\a}, \quad \theta \in \mathbb{R}. 
    \end{equation}
    Moreover, for every $0<\alpha \leq 2$, $\Lambda_\a$ is self-decomposable; see e.g. \cite[page 192]{YYY09}. A direct calculation gives
    \begin{equation}
        J_{\Lambda_\a }(z) = \frac{1}{1+|z|^\alpha}.
    \end{equation}
    But $ J_{\Lambda_\a }(z) \in L^1(\mathbb{R})$ if $1 < \a \leq 2$ and $ J_{\Lambda_\a }(z) \notin L^1(\mathbb{R})$ if $0 < \a \leq 1$.
\end{remark}

\section{Proof of Theorem \ref{Theorem non-SD of alpha Cauchy}}  
\label{sec proof non-SD}
The particular case $\mathcal{C}_2$ is the well-known standard Cauchy variable, which is infinitely divisible and self-decomposable. Recently, Wang \cite{Wan26} proved that $\mathcal{C}_\alpha$ is not infinitely divisible if $\alpha > 2$, thus $\mathcal{C}_\alpha$ is not self-decomposable if $\alpha > 2$. In the following, we focus on the case $1 < \alpha <2$. 

Let 
$$\phi_\a(z) := \mathbb{E}[e^{iz \mathcal{C}_\a }] \quad \text{and} \quad  \Psi_\a(z) := -\log \phi_\a (z). $$
By \eqref{eq phi' sym at 0}, we have
\begin{equation}
    z \Psi_\a'(z) \sim   k_\a \Gamma(2-\alpha) |z|^{\alpha -1} ,  \quad z \rightarrow 0. 
\end{equation}
Because $1< \alpha < 2$, $\lim_{z \rightarrow  0} z\Psi_\a'(z) = 0$. \\
By \eqref{eq phi infty} and \eqref{eq phi' infty}, we have 
\begin{equation}
  \lim_{z \rightarrow  \infty}  z \Psi_\a'(z)    = 1 + \alpha. 
\end{equation}

Then we define 
\begin{equation}
    \label{eq def J_alpha}
    J_\a(z) := 1 - \frac{z \Psi_\a'(z)}{1+\a}.
\end{equation}
\begin{center}
    Claim 1: $\int_\mathbb{R}|J_\a(z)|dz < \infty,  \quad$ and $\quad  \int_\mathbb{R} J_\a (z)dz < 0.$
\end{center}
Theorem \ref{theorem new criterion SD} and Claim 1 imply that 
$\mathcal{C}_\a$ is not self-decomposable. 

It suffices to prove Claim 1.

\subsection{Proof of Claim 1} 
The function $z \mapsto J_\a(z)$ is even, since $z \mapsto z \Psi_\a'(z)$ is even. \\
We want to prove 
$\int_0^\infty|J_\a(z)|dz < \infty,$ and $  \int_0^\infty J_\a (z)dz < 0.$
By \eqref{eq phi infty} and \eqref{eq phi' infty}, we have
\begin{equation}
\label{eq sym at infty}
    J_\a(z) \sim  - \frac{r_\alpha}{1+\alpha} z ^{-\alpha} , \quad z \rightarrow +\infty, 
\end{equation}
where $r_\alpha = - 2  \cos(\pi \alpha/2)  \alpha \Gamma(1+2\alpha)/\Gamma(1+\alpha) > 0.$ 
Therefore, $\int_0^\infty|J_\a(z)|dz < \infty.$  

\textbf{It remains to prove $\int_0^\infty J_\a (z)dz < 0.$ }

\begin{align}
    \int_0^\infty J_\a (z)dz = &  \int_0^\infty (1 - \frac{z\Psi'_\a (z)}{1 + \a}) dz \\
    = &  \int_0^\infty \frac{z}{1 + \a}(\frac{1 + \a}{z} - \Psi'_\a (z)) dz   \\
    = &  \int_0^\infty \frac{z}{1 + \a} d[(1 + \a)\log z - \Psi_\a (z) - \log (k_\alpha \Gamma(1+\alpha))]  \\
    = &  \int_0^\infty \frac{z}{1 + \a} d \left[\log \frac{z^{1+\a} \phi_\a(z)}{k_\alpha \Gamma(1+\alpha)} \right]   \\
    = &  \lim_{z\rightarrow \infty}  \frac{z}{1 + \a} \log \frac{z^{1+\a} \phi_\a(z)}{k_\alpha \Gamma(1+\alpha)}  - \lim_{z\rightarrow 0}  \frac{z}{1 + \a} \log \frac{z^{1+\a} \phi_\a(z)}{k_\alpha \Gamma(1+\alpha)}\\
    \label{eq part}
    & - \int_0^\infty  \left[\log \frac{z^{1+\a} \phi_\a(z)}{k_\alpha \Gamma(1+\alpha)} \right]  \frac{dz}{1 + \a}. 
\end{align}
\eqref{eq phi infty} implies that 
\begin{equation}\label{eq log-expansion}
\log\frac{z^{1+\alpha}\phi_\alpha(z)}
 {k_\alpha\Gamma(1+\alpha)}
 =-2\cos(\pi\alpha/2)
 \frac{\Gamma(1+2\alpha)}{\Gamma(1+\alpha)}z^{-\alpha}
 +O(z^{-2\alpha}), \quad z \rightarrow \infty.
\end{equation}
Because $1< \alpha <2$, we have
\begin{equation}
     \lim_{z\rightarrow \infty}  \frac{z}{1 + \a} \log \frac{z^{1+\a} \phi_\a(z)}{k_\alpha \Gamma(1+\alpha)} = 0. 
\end{equation}
Because $\phi_\a(0 ) = 1$, we have 
\begin{equation}
     \lim_{z\rightarrow 0}  \frac{z}{1 + \a} \log \frac{z^{1+\a} \phi_\a(z)}{k_\alpha \Gamma(1+\alpha)} = 0. 
\end{equation}
Therefore, by \eqref{eq part}, the inequality $\int_0^\infty J_\a (z)dz < 0$ is equivalent to 
\begin{align}
  0 > & - \int_0^\infty  \log \frac{z^{1+\a} \phi_\a(z)}{k_\alpha \Gamma(1+\alpha)}dz \\
  = &  - \int_0^\infty \log \left[ \frac{z^{1+\a}}{\Gamma(1+\alpha)} \int_0^\infty  e^{-zy}  \frac{ y^\alpha}{  y^{2\alpha} + 2\cos(\pi \alpha/2)y^\alpha + 1}dy \right] dz  \\
    = &  -\int_0^\infty \log \left[ \frac{1}{\Gamma(1+\alpha)} \int_0^\infty  e^{-u}  \frac{ u^\alpha}{  (u/z)^{2\alpha} + 2\cos(\pi \alpha/2)(u/z)^\alpha + 1}du \right] dz \\
    \label{eq bound}
    = & -\int_0^\infty \log \mathbb{E}\left[ H_\alpha(\frac{\G_{1+\alpha}}{z})\right] dz, 
\end{align}
where $$   H_\alpha(y)
    :=\frac{1}{y^{2\alpha}+2\cos(\pi \alpha/2)y^\alpha+1},
    \qquad y>0,$$
and $\G_{c}$ denotes the Gamma random variable, whose density is $ \frac{1}{\Gamma(c)}x^{c-1}e^{-x}\mathbf{1}_{(0, \infty)}(x). $

For every fixed $z> 0$, $H_\alpha(\frac{\G_{1+\alpha}}{z})$ is an integrable random variable and $-\log$ is a convex function, we use Jensen's inequality to obtain 
\begin{equation}
    -\log \mathbb{E}\left[ H_\alpha(\frac{\G_{1+\alpha}}{z})\right] \leq   \mathbb{E}\left[ -\log H_\alpha(\frac{\G_{1+\alpha}}{z})\right].
\end{equation}
Note that the convex function $-\log$ is strictly convex, and the random variable $H_\alpha(\frac{\G_{1+\alpha}}{z})$ is not a constant, the above equality never holds, i.e. we have for every fixed $z > 0$, 
\begin{equation}\label{eq bound2}
    -\log \mathbb{E}\left[ H_\alpha(\frac{\G_{1+\alpha}}{z})\right] <   \mathbb{E}\left[ -\log H_\alpha(\frac{\G_{1+\alpha}}{z})\right].
\end{equation}
We next prove 
\begin{equation}
    \int_0^\infty  \mathbb{E}\left[ -\log H_\alpha(\frac{\G_{1+\alpha}}{z})\right] dz  = 0; 
\end{equation}
if this equality holds, combining \eqref{eq bound2}  and \eqref{eq bound}, we finish the proof. 
\begin{align}
   &  \int_0^\infty  \mathbb{E}\left[ -\log H_\alpha(\frac{\G_{1+\alpha}}{z})\right] dz \\
 = &   \int_0^\infty \left[ \int_0^\infty  -\log H_\alpha(\frac{u}{z}) \frac{1}{\Gamma(1+\alpha) }u^\alpha e^{-u} du\right] dz\\
= &    \int_0^\infty \left[ \int_0^\infty  \log ((u/z)^{2\alpha} + 2\cos(\pi \alpha/2)(u/z)^\alpha + 1) \frac{1}{\Gamma(1+\alpha) }u^\alpha e^{-u}du \right] dz\\
= &   \int_0^\infty \left[ \int_0^\infty  \log (y^{2\alpha} + 2\cos(\pi \alpha/2)y^\alpha + 1) \frac{1}{\Gamma(1+\alpha) }(yz)^\alpha e^{-yz} zdy \right] dz \\
\label{eq using Fubini}
= & \int_0^\infty \left[ \int_0^\infty  \log (y^{2\alpha} + 2\cos(\pi \alpha/2)y^\alpha + 1) \frac{1}{\Gamma(1+\alpha) }(yz)^\alpha e^{-yz} zdz \right] dy  \\
\label{eq last step}
= &  \int_0^\infty  \log (y^{2\alpha} + 2\cos(\pi \alpha/2)y^\alpha + 1) \frac{1+\alpha}{y^2}  dy. 
\end{align}
The equality \eqref{eq using Fubini} follows from Fubini's Theorem, whose application is guaranteed by
\begin{align}
    &  \int_0^\infty \left[ \int_0^\infty | \log (y^{2\alpha} + 2\cos(\pi \alpha/2)y^\alpha + 1)| \frac{1}{\Gamma(1+\alpha) }(yz)^\alpha e^{-yz} zdy \right] dz  \\
    \label{eq using Tonelli}
    = &  \int_0^\infty \left[ \int_0^\infty | \log (y^{2\alpha} + 2\cos(\pi \alpha/2)y^\alpha + 1)| \frac{1}{\Gamma(1+\alpha) }(yz)^\alpha e^{-yz} zdz \right] dy \\
    = &   \int_0^\infty  |\log (y^{2\alpha} + 2\cos(\pi \alpha/2)y^\alpha + 1)| \frac{1+\alpha}{y^2}  dy \\
    \label{eq absolutely integrable}
    < & \infty. 
\end{align}
The equality \eqref{eq using Tonelli} follows from Tonelli's Theorem. The inequality \eqref{eq absolutely integrable} holds because the integrand is of order $y^{\alpha -2}$ when $y \downarrow 0$ and of order $(\log y)/y^2$ when $y \rightarrow +\infty$.

Finally, by \eqref{eq last step} we wish to prove 
\begin{equation}
     \int_0^\infty  \log (y^{2\alpha} + 2\cos(\pi \alpha/2)y^\alpha + 1) \frac{dy}{y^2}   =  0.   
\end{equation}
By direct computation, 
\begin{align}
   &  \int_0^\infty  \log (y^{2\alpha} + 2\cos(\pi \alpha/2)y^\alpha + 1) \frac{dy}{y^2} \\
   = & \int_0^\infty  \log (y^{2\alpha} + 2\cos(\pi \alpha/2)y^\alpha + 1) (-dy^{-1} ) \\
   = & \lim_{y \rightarrow 0} \frac{ \log (y^{2\alpha} + 2\cos(\pi \alpha/2)y^\alpha + 1)}{y} - \lim_{y \rightarrow \infty} \frac{ \log (y^{2\alpha} + 2\cos(\pi \alpha/2)y^\alpha + 1)}{y} \\
   & + \int_0^\infty \frac{1}{y} d\log (y^{2\alpha} + 2\cos(\pi \alpha/2)y^\alpha + 1) \\
   = & \alpha \int_0^\infty  \frac{ 2y^{2\alpha-2} + 2\cos(\pi \alpha/2) y^{\alpha-2}}{y^{2\alpha} + 2\cos(\pi \alpha/2)y^\alpha + 1}dy = \int_0^\infty  \frac{ 2 u + 2\cos(\pi \alpha/2) }{u^{2} + 2\cos(\pi \alpha/2)u + 1} u^{-1/\alpha} du \\
   = &  \int_0^\infty \frac{u^{-1/\alpha} e^{i\pi \alpha/2} }{1 + e^{i\pi \alpha/2} u }+ \frac{u^{-1/\alpha} e^{- i\pi \alpha/2} }{1 + e^{-i\pi \alpha/2} u }   du   \\
   = &  e^{i\pi \alpha/2} \int_0^\infty \frac{u^{-1/\alpha}  }{1 + e^{i\pi \alpha/2} u }  du  + e^{-i\pi \alpha/2} \int_0^\infty \frac{u^{-1/\alpha}  }{1 + e^{-i\pi \alpha/2} u }  du \\
   \label{eq residue}
   =& e^{i\pi/2} \int_0^\infty \frac{u^{-1/\alpha}  }{1 +  u }  du + e^{-i\pi/2} \int_0^\infty \frac{u^{-1/\alpha}  }{1 +  u }  du \\
   =&  0. 
\end{align}
The equality \eqref{eq residue} follows from using Residue theorem for those two integral separately. This completes the proof.

\section{Proof of lemma \ref{lemma asymptotic estimates alpha Cauchy}}
\label{sec proof lemma}
We first observe that
\begin{equation}
 \frac{1}{1+2\cos(\pi\alpha/2) y^\alpha+y^{2\alpha}}
 =1-2\cos(\pi\alpha/2) y^\alpha+O(y^{2\alpha}),
 \qquad y\downarrow0.
\end{equation}
Apply Watson's lemma (see (2.3.7) and (2.3.8) in \cite{NIST}) separately to the integral representations of
$\phi_\alpha$ and $\phi_\alpha'$, we have \eqref{eq phi infty} and \eqref{eq phi' infty}. 
By Hardy–Littlewood Tauberian theorem, we have
\eqref{eq phi' sym at 0}. \\
For the reader's convenience, we rewrite Hardy–Littlewood Tauberian theorem (see, e.g. \cite[formula (3.2)]{MS13}) here. 
Let $F:[0,\infty)\to\mathbb{R}$ be a nondecreasing right-continuous function. And suppose that its Laplace-Stieltjes transforms
$$\omega(s) = \int_0^\infty e^{-st}\,dF(t)$$
exists for some $s > 0$. 
The theorem states that 
$$\omega(s)\sim C s^{-\rho},\quad\rm{as\ }s\downarrow 0 \quad \Longleftrightarrow \quad  F(t)\sim \frac{C}{\Gamma(\rho+1)}t^\rho, \ \text{as}\ t\rightarrow +\infty$$
where $C > 0$ and $\rho \geq 0.$

\section{Proof of Corollary \ref{Corollary}}
\label{sec proof corollary}
Recall that the $\alpha$-Cauchy variables were introduced to study the first hitting times of points for one-dimensional symmetric stable L\'evy processes. 
Let $X_\alpha = (X_\alpha(t): t\geq 0)$ be the symmetric stable L\'evy process of index $\alpha$ starting from zero and \(\hat{X}_\alpha = (\hat{X}_\alpha(t) : t \geq 0)\) be an independent copy of \(X_\alpha\).
Let $T_{\{a\}}(X_\alpha)$ denote the first hitting time of point $a$ for $X_\alpha$. We rewrite formula (5.12) in \cite{YYY09} here
 \begin{equation}
    \label{relation1}
\hat{X}_{\alpha}(T_{\{a\}}(X_{\alpha})) \overset{\text{law}}{=} |a| \mathcal{C}_{\alpha}. 
 \end{equation}
It is known that if the subordinand $\{X_t\}$ is strictly stable and the subordinator $\{Z_t\}$ is self-decomposable, and they are independent, then the subordinated process ${X_{Z_t}}$ is self-decomposable, see e.g. Theorem 5.1 in \cite{AS19}. 
Using this subordination property, the formula \eqref{relation1} and Theorem \ref{Theorem non-SD of alpha Cauchy}, we obtain Corollary \ref{Corollary}.

\begin{remark}
We emphasize, in the proof of Corollary \ref{Corollary}, that the subordinand is strictly stable. 
    There exists a stable process $(X(t): t\geq 0)$ and a self-decomposable subordinator $(Z(t): t\geq 0)$, such that the corresponding subordinated process ${X_{Z_t}}$ is not self-decomposable; see \cite{Koz05}. 
\end{remark}

\bigskip

\noindent
\textbf{Acknowledgements.} We thank Thomas Simon for reading the first draft of this paper and for providing the article \cite{Koz05}.

\bigskip

\noindent
\textbf{Declaration of AI Use.}
The authors used ChatGPT (OpenAI) for assistance with language editing and for discussing and checking mathematical arguments during the preparation of this manuscript. All mathematical statements, arguments, and proofs included in the final version were independently verified by the authors, who take full responsibility for the content of the article.


\bibliographystyle{plain} 
\bibliography{ID}
\end{document}